\documentstyle[art10]{article}
\title{Projective normality of special scrolls.}
\author{Luis Fuentes Garc\'{\i}a\thanks{Supported by an F.P.U.
fellowship of Spanish Government}
 \and{Manuel Pedreira P\'erez}
 }
\date{}

\newtheorem{teo}{Theorem}[section]
\newtheorem{defin}[teo]{Definition}
\newtheorem{prop}[teo]{Proposition}
\newtheorem{cor}[teo]{Corollary}
\newtheorem{lemma}[teo]{Lemma}

\newtheorem{rem}[teo]{Remark}

\font\euf=eufm10 at 12pt

\def\g2{\pi}

\def\e{\mbox{\euf e}}

\def\b{\mbox{\euf b}}

\def\a{\mbox{\euf a}}

\def\k{{\cal K}}
\def\K{{\cal K}}

\def\P{{\bf P}}

\newcommand\E{{\cal E}}
\newcommand\Te{{\cal O}}

\def\qed{\hspace{\fill}$\rule{2mm}{2mm}$}
\newcommand\lrw{\longrightarrow}
\newcommand\rw{\rightarrow}

\newcommand{\Sym}{\mathop{\rm Sym}\nolimits}

\def\ov{\overline}

\begin{document}
\maketitle

{\footnotesize{\bf Authors' address:} Departamento de Algebra, Universidad de Santiago
de Compostela. $15706$ Santiago de Compostela. Galicia. Spain. e-mail: {\tt
pedreira@zmat.usc.es}; \\ {\tt luisfg@usc.es}\\
{\bf Abstract:}  We study the projective normality of a linearly normal special scroll $R$ of degree
$d$ and speciality $i$ over a smooth curve $X$ of genus $g$. If $d\geq 4g-2i+1$, we prove that the
projective normality of the scroll is equivalent to the projective normality of its directrix curve
of minimum degree.}

\vspace{0.1cm}

{\Large\bf Introduction.} 

Let $R\subset \P^N$ be a scroll, that is, a surface with a one-dimensional family of lines. We say that
$R$ is linearly normal when $R$ is not the projection of any other scroll with the same degree and genus.
It is known that there are a unique geometrically ruled surface $\pi:S=\P(\E)\lrw X$ and a unique base
point free unisecant linear system $|H|$ on $S$ such that $R$ is the birational image of $S=\P(\E)$ by
the regular map defined by the complete linear system $|H|$. We define the speciality of $R$ as
$i=h^1(\Te_S(H))$. The scroll is special when $i>0$.

In this paper we study the projective normality of a special scroll $R$. We say that $R$ is projectively
normal, or equivalently, that the line bundle $\Te_S(H)$ on $S$ is normally generated, if and only if the
natural maps $Sym_k(H^0(\Te_S(H)))\lrw H^0(\Te_S(kH))$ are surjective for any
$k\geq 1$. We say that the cokernel of these maps are the speciality of
$R$ respect to hypersurfaces of degree $k$. 

The question of the normal generation of a line bundle ${\cal L}$ on a smooth curve $X$ has been studied
with detail.  The classical result of Noether says that the canonical bundle is normally generated;
the Theorem of Castelnuovo states that any line bundle of degree at least $2g+1$ is normally generated.
More recently, H. Lange and G. Martens in \cite{lange} and M. Green and R. Lazarsfeld in \cite{green}
gave bounds, depending of the Clifford index of $X$, over the degree of a line bundle to be normally
generated. 

There are not complete generalizations of these theorems to higer dimension. If $V$ is a
projective variety of dimension $n$, most of the results study the normal generation of adjoint
bundles, that is, line bundles $\omega_V\otimes {\cal L}$ where ${\cal L}$ is an ample line bundle.

M. Andreatta and J. Sommese in \cite{andreatta} and L.Ein and R.Lazarsfled in  \cite{ein} proved that
$\omega_V\otimes {\cal L}^{\otimes n}\otimes {\cal F}$ is normally generated if it is very ample,
${\cal L}$ is very ample and ${\cal F}$ is numerically effective. D. Butler in \cite{butler} has
stabilished that if
$V$ is a ruled variety of rank $n$ over a smooth curve $X$, then
$\omega_V\otimes {\cal L}^{\otimes 2n+1}$  is normally generated if ${\cal L}$ is ample.

However, the canonical divisor of a ruled surface $S$ is $\K_S\sim -2X_0+(\k+\e)f$, so in general, it is
not easy to write an unisecant linear system like $\K_S+2A+B$ with $|A|$ very ample and $|B|$ numerically
effective. Furhtermore, an ample divisor on $S$ must be at least unisecant so the 
result of D. Butler applies at least to $3$-secant linear systems on $S$. Then, we
will use more particular strategies to the case of ruled surfaces.  We try to reduce the study of the
normal generation of a line bundle
$\Te_S(H)$ to the normal generation of line bundles on the base curve $X$. We will prove that:

{\em The projective normality of a special scroll of genus $g$, speciality $i$ and degree $d\geq
4g-2i+1$ is equivalent to the projective normality of its directrix curve of minimum degree.}

We refer to \cite{fuentes} for a systematic
development of the projective theory of scrolls and ruled surfaces that we will use in this paper and to
\cite{fuentes2} to study the special scrolls. Anyway, in the first section we recall some
basic facts that we will use along the paper.

\section{Preliminars.}

A {\it geometrically ruled surface}, or simply a {\it
ruled surface}, will be a $\P^1$-bundle over a smooth curve $X$ of genus $g>0$. It will be
denoted by $\pi: S=\P(\E_0)\lrw X$. We will suppose that
$\E_0$ is a normalized sheaf and $X_0$ is the section of minimum self-intersection
that corresponds to the surjection $\E_0\lrw \Te_X(\e)\lrw 0$, $\bigwedge^2\E\cong
\Te_X(\e)$ (see \cite{hartshorne},V, \S 2 and \cite{fuentes}).

If $|H|=|X_0+\b f|$ is a base-point-free linear system on a ruled surface $S$, $|H|$ defines a regular
map $\phi_H:S\lrw \P(H^o(\Te_S(H)^\vee)$. The image of $S$ is a scroll $R$. If $\phi_H$ is a birational
map we say that
$S$ and $H$ are the ruled surface and the linear system associated to the scroll $R$. We denote de image
of a curve $C\subset S$ by $\ov{C}\subset R$. The curve $\ov{X_0}$ is the curve of minimum degree of $R$.
It is embedded by the linear system $|\b+\e|$ on $X$.

We will use the well known fact that, if $m\geq 0$ and $i\geq 0$:
$$
h^i(\Te_S(mX_0+\b f))\leq\sum\limits^m_{k=0} h^i(\Te_X(\b + k\e))
$$
Furthermore, if $\b$ is nonspecial $$h^i(\Te_S(X_0+\b
f))=h^i(\Te_X(\b))+h^i(\Te_X(\b+\e))$$ (see \cite{fuentes}).

\begin{defin}
Let $V$ be a projective variety. Let ${\cal F}_i,i=1,\ldots ,s$, coherent sheaves on $V$. We
call
$s({\cal F}_1,\ldots,{\cal F}_s)$ the cokernel of the map:
$$
H^0({\cal F}_1)\otimes \ldots \otimes H^0({\cal F}_s)\lrw H^0({\cal F}_1\otimes \ldots \otimes {\cal
F}_s)
$$
If ${\cal F}_i$ are invertible sheaves $\Te_V(D_i)$ where $D_i$ are divisors on $V$, we will write
$s(D_1,\ldots,D_s)$.
\end{defin}

\begin{lemma}\label{lematecnico}
If $s({\cal F}_1,{\cal F}_2)=0$, then:
$$
s({\cal F}_1,{\cal F}_2,{\cal F}_3,\ldots,{\cal F}_s)=s({\cal F}_1\otimes {\cal
F}_2,{\cal F}_3,\ldots,{\cal F}_s)
$$
\end{lemma}
{\bf Proof:}  It is sufficient note that $s({\cal F}_1,{\cal F}_2,{\cal F}_3,\ldots,{\cal F}_s)$ is the
cokernel of the composition:
$$
H^0({\cal F}_1)\otimes \ldots \otimes H^0({\cal F}_s)\rw
H^0({\cal F}_1\otimes {\cal F}_2)\otimes H^0({\cal F}_3)\otimes \ldots \otimes H^0({\cal F}_s)
\rw H^0({\cal F}_1\otimes \ldots \otimes {\cal F}_s)
$$ \qed

\begin{defin}

Let $V$ be a projective variety and let $|H|$ be a complete unisecant base-point-free linear system
defining a birational map:
$$
\phi_H:V\lrw \ov{V}\subset \P^N
$$
We say that $(V,\Te_V(H))$ is projectively normal or $\Te_V(H)$ normally generated or $\ov{V}$
projectively normal if and only if the natural maps:
$$
\Sym_k(H^0(\Te_V(H)))\lrw H^0(\Te_V(kH))
$$
are surjective for all $k\geq 1$.

\end{defin}

\begin{rem}
Thus, $(V,\Te_V(H))$ is projectively normal if and only if $s(H,\stackrel{k}{\ldots},H)=0$, for all
$k\geq 2$. Moreover, $dim(s(H,\stackrel{k}{\ldots},H))$ is the speciality of $\ov{V}$ respect to
hypersurfaces of degree $k$. 
\end{rem}

\begin{lemma}[Green]\label{green}
Let $\a$,$\b$ be effective divisors on a smooth curve $X$. Let $\b$ be base-point-free. If
$h^1(\Te_X(\a-\b))\leq h^0(\Te_X(\b))-2$ then $s(\a,\b)=0$.
\end{lemma}
{\bf Proof:} This is a particular case of $H^o$-Lemma in \cite{lgreen}.

\section{Projective normality of special scrolls.}

\begin{prop}\label{nproyectivascroll}
Let $S$ be a ruled surface and $H\sim X_0+\b f$ an unisecant linear system, such that $\b$ is a
nonspecial divisor and $\b+\e$ is effective. Let $k\geq 2$. If
$$
s(\b+\e,\stackrel{i}{\ldots},\b+\e,\b,\stackrel{k-i}{\ldots},\b)=0 \mbox{ for all $i$, with $0\leq i\leq
k-1$}
$$
then:
$$
s(H,\stackrel{k}{\ldots},H)\cong s(\b+\e,\stackrel{k}{\ldots},\b+\e)
$$
\end{prop}
{\bf Proof:} Let $i,j$ be non negative integers, with $i+j>0$. We will denote by
$W_{i,j}$ the tensor product:
$$
H^0(\Te_S(H))\otimes \stackrel{i}{\ldots}\otimes H^0(\Te_S(H))\otimes
	   H^0(\Te_S(H-X_0))\otimes \stackrel{j}{\ldots}\otimes H^0(\Te_S(H-X_0))	
$$
With this notation it verifies that:
\begin{enumerate}

\item $W_{i,j}\otimes H^0(\Te_S(H))=W_{i+1,j}$.

\item $W_{i,j}\otimes H^0(\Te_S(H-X_0))=W_{i,j+1}$

\item We have an exact sequence:
$$
W_{i,j}\rw H^0(\Te_S((i+j)H-jX_0)))\rw
s(H,\stackrel{i}{\ldots},H,H-X_0,\stackrel{j}{\ldots},H-X_0)\rw 0
$$

\item The map $H^0(\Te_S(H))\lrw H^0(\Te_{X_0}(H))$ is a surjection, because $\b$ is nonspecial.
Moreover, we know that $H^0(\Te_S(S-X_0))\cong H^0(\Te_X(\b))$ and $H^0(\Te_{X_0}(H))\cong 
H^0(\Te_X(\b+\e))$. From this, we have the following exact sequence:
$$
\begin{array}{l}
{W_{i-1,j}\otimes H^0(\Te_{X_0}(H))\lrw H^0(\Te_X((i+j)\b+i\e))\lrw}\\
{\lrw s(\b+\e,\stackrel{i}{\ldots},\b+\e,\b,\stackrel{j}{\ldots},\b)\rw 0}\\
\end{array}
$$
\end{enumerate}

We will prove that $s(H,\stackrel{i}{\ldots},H,H-X_0,\stackrel{k-i}{\ldots},H-X_0)=0$
if $0\leq i\leq k-1$.

We proceed by induction on $i$:

\begin{enumerate}

\item If $i=0$, because $H^0(\Te_S(H-X_0))\cong H^0(\Te_X(\b))$ we inmediately obtain that:
$$
s(H-X_0,\stackrel{k}{\ldots},H-X_0)=s(\b,\stackrel{k}{\ldots},\b)
$$
and this is zero by hypotesis.

\item Now, suppose that the assertion is true for $i-1$, with $1\leq i\leq k-1$.
Consider the exact sequence:
$$
0\lrw H^0(\Te_S(H-X_0))\lrw H^0(\Te_S(H))\lrw H^0(\Te_{X_0}(H))\lrw 0
$$
Taking the tensor product with  $W_{i-1,k-i}$ we obtain the following conmutative diagram:
$$
\setlength{\unitlength}{5mm}
\begin{picture}(22,9)

\put(4,1){\makebox(0,0){$0$}}
\put(4,3){\makebox(0,0){$W_{i-1,k-i}\otimes H^0(\Te_{X_0}(H))$}}
\put(4,5){\makebox(0,0){$W_{i,k-i}$}}
\put(4,7){\makebox(0,0){$W_{i-1,k-(i-1)}$}}
\put(4,9){\makebox(0,0){$0$}}

\put(17,1){\makebox(0,0){$H^1(\Te_S(kH-(k-(i-1))X_0))$}}
\put(17,3){\makebox(0,0){$H^0(\Te_X(k\b+i\e))$}}
\put(17,5){\makebox(0,0){$H^0(\Te_S(kH-(k-i)X_0))$}}
\put(17,7){\makebox(0,0){$H^0(\Te_S(kH-(k-(i-1))X_0))$}}
\put(17,9){\makebox(0,0){$0$}}

\put(4,2.5){\vector(0,-1){1}}
\put(4,4.5){\vector(0,-1){1}}
\put(4,6.5){\vector(0,-1){1}}
\put(4,8.5){\vector(0,-1){1}}

\put(17,2.5){\vector(0,-1){1}}
\put(17,4.5){\vector(0,-1){1}}
\put(17,6.5){\vector(0,-1){1}}
\put(17,8.5){\vector(0,-1){1}}

\put(8.5,3){\vector(1,0){3}}
\put(8.5,5){\vector(1,0){3}}
\put(8.5,7){\vector(1,0){3}}

\put(10,3.5){\makebox(0,0){$\alpha_{3i}$}}
\put(10,5.5){\makebox(0,0){$\alpha_{2i}$}}
\put(10,7.5){\makebox(0,0){$\alpha_{1i}$}}

\end{picture}
$$
where
$$
\begin{array}{l}
{coker(\alpha_{1i})=s(H,\stackrel{i-1}{\ldots},H,H-X_0,\stackrel{k-(i-1)}{\ldots},H-X_0)}\\
{coker(\alpha_{2i})=s(H,\stackrel{i}{\ldots},H,H-X_0,\stackrel{k-i}{\ldots},H-X_0)}\\
{coker(\alpha_{3i})=s(\b+\e,\stackrel{i}{\ldots},\b+\e,\b,\stackrel{k-i}{\ldots},\b)}\\
\end{array}
$$
Moreover,
$$
h^1(\Te_S(kH-(k-(i-1))X_0))\leq \sum\limits^{i-1}_{l=0} h^1(\Te_X((k-l)\b + l(\b+\e)))
$$
Because $\b$ is nonspecial and $\b+\e$ is effective, if $i\leq k$ and $0\leq l\leq i-1$ then  $(k-l)\b +
l(\b+\e)$ is nonspecial. Therefore, $h^1(\Te_S(kH-(k-(i-1))X_0)))=0$.

By induction hypothesis, $s(H,\stackrel{i-1}{\ldots},H,H-X_0,\stackrel{k-(i-1)}{\ldots},H-X_0)=0$ and by
the hypothesis of the theorem
$s(\b+\e,\stackrel{i}{\ldots},\b+\e,\b,\stackrel{k-i}{\ldots},\b)=0$ when
$i\leq k-1$. Thus, since $\alpha_{1i}$ and $\alpha_{3i}$ are surjections, we deduce that
$\alpha_{2i}$ is a surjection, and then
$s(H,\stackrel{i}{\ldots},H,H-X_0,\stackrel{k-i}{\ldots},H-X_0)=0$.

\end{enumerate}

Finally, from the above diagram for $i=k$, we see that $\alpha_{1k}$ is a surjection, because
$s(H,\stackrel{k-1}{\ldots},H,H-X_0)=0$. Therefore, $coker(\alpha_{2k})\cong coker(\alpha_{3k})$ and,
$$
s(H,\stackrel{k}{\ldots},H)\cong s(\b+\e,\stackrel{k}{\ldots},\b+\e)
$$
\qed

\begin{cor}
Let $S$ be a geometrically ruled surface. Let $|H|=|X_0+\b f|$ an unisecant complete linear system
on $S$, such that $\b$ is a nonspecial divisor and $\b+\e$ is effective. If
$$
s(\b+\e,\stackrel{i}{\ldots},\b+\e,\b,\stackrel{k-i}{\ldots},\b)=0 \mbox{ for all $i$, with $0\leq i\leq
k$}
$$
for all $k\geq 2$ then  $(S,\Te_S(H))$ is projectively normal.
\end{cor}

\begin{prop}
Let  $R\subset \P^N$ be a special linearly normal scroll of genus $g$, degree $d$ and speciality $i$.
Let $S$ and $|H|=|X_0+\b f|$ be the ruled surface and the linear systema asocciated to $R$.

If $\b$ is nonspecial then the curve of minimum degree of
$R$ is special, linearly normal and its speciality coincides with the specialty of $R$. Moreover, if
$deg(\b)\geq 2g-1$ then this is the unique special curve of the scroll.
\end{prop}
{\bf Proof:} If $\b$ is nonspecial we know that:
$$
h^1(\Te_S(H))=h^1(\Te_X(\b))+h^1(\Te_X(\b+\e))
$$
Thus, the speciality of the curve of minimun degree $\ov{X_0}$ is $h^1(\Te_X(\b+\e))=h^1(\Te_S(H))=i$.
Furthermore, $\ov{X_0}$ is linearly normal because we have a surjective map: 
$$
H^0(\Te_S(H))\lrw H^0(\Te_{X_0}(H))\lrw H^1(\Te_S(\b))=0
$$

If $deg(\b)\geq 2g-1$, the degree of the scroll satisfies $d\geq 2g-1+deg(\ov{X_0})$. Any other
directrix curve of $R$ different from
$\ov{X_0}$ has degree $d'\geq d-deg(\ov{X_0})=2g-1> 2g-2$, so it is nonspecial. \qed

\begin{teo}\label{directrizespecialminima1}
Let  $R\subset \P^N$ be a special linearly normal scroll of genus $g$, degree $d$ and speciality $i$. 
Suppose that $d\geq 4g-2i+1$. Then:
\begin{enumerate}

\item $R$ has an unique special directrix curve $\ov{X_0}$. Moreover,
$\ov{X_0}$ is the curve of minimum degree, it is linearly normal and it has the speciality of $R$.

\item $R$ and $\ov{X_0}$ have the same specialty respect to hypersurfaces of degree $m$. In particular
the scroll is projectively normal if and only if the curve of minimun degree if projectively normal.

\end{enumerate}

\end{teo}
{\bf Proof:} Let  $S$ be the ruled surface and  $|H|=|X_0+\b f|$ the linear system corresponding to the
scroll $R$.

\begin{enumerate}

\item Since  $R$ is special, it has a special directrix curve (see \cite{fuentes2}) so the curve
$\ov{X_0}$ of minimun degree of the scroll  verifies $deg(\b+\e)\leq 2g-2$.

Suppose that $d\geq 4g-2i+1$. We know that:
$$
d-2g+2+i=h^0(\Te_S(H))\leq h^0(\Te_X(\b+\e))+h^0(\Te_X(\b))
$$

If $deg(\b)\leq 2g$, then we can apply Clifford Theorem  (\cite{hartshorne}, page 343) to the
divisors
$\b$ and
$\b+\e$. From the above inequiality we obtain:
$$
d-2g+2+i\leq \frac{deg(\b+\e)}{2}+1+\frac{deg(\b)}{2}+1=\frac{d}{2}+2
$$
and then  $d\leq 4g-2i$. 

Thus, we see that $deg(\b)\geq 2g+1$ and now we can apply the Proposition above.

\item We will apply the Proposition \ref{nproyectivascroll}. We have to proof that:
$$
s(\b+\e,\stackrel{i}{\ldots},\b+\e,\b,\stackrel{k-i}{\ldots},\b)=0 \mbox{ with $0\leq i\leq k-1$}
$$
for all $k\geq 2$.

If $\b+\e\sim 0$, then the scroll is a cone and it verifies the condition trivialy. Thus, we will
suppose that $\b+\e\not\sim 0$.

Since $deg(\b)\geq 2g+1$, we can use the Theorem of Castelnuovo to see that $\Te_X(\b)$ is normaly
generated and then $s(\b,\stackrel{k}{\ldots},\b)=0$ for any $k\geq 2$. From this, applying
Lemma \ref{lematecnico}:
$$
s(\b,\stackrel{i}{\ldots},\b,\b+\e,\stackrel{k-i}{\ldots},\b+\e)=s(i\b,\b+\e,\stackrel{k-i}{\ldots},\b+\e)
$$
when $1\leq i\leq k$.

Let us see that $s(\b,\b+\e)=0$. We apply  the Lemma
\ref{green}. We need that:
$$
h^1(\Te_X(\b-(\b+\e)))< h^0(\Te_X(\b+\e))-1
$$
Note that $e=-deg(\e)\geq 0$ because  $deg(\b)\geq 2g+1$ and $deg(\b+\e)\leq 2g-2$.
We distinguish to cases:
\begin{enumerate}

\item If $-\e$ is nonspecial, it is sufficient to check that $h^0(\Te_X(\b+\e))>
1$. But $|\b+\e|$ is base point free and $\b+\e\not\sim 0$, so $h^0(\Te_X(\b+\e))> 1$. 

\item If $-\e$ is special, we apply Clifford Theorem  (\cite{hartshorne}, page 343):
$$
|-\e|\leq \frac{1}{2} e \mbox{,  or equivalently,  } h^1(\Te_X(-\e))\leq g-\frac{e}{2}
$$
But,
$$
h^0(\Te_X(\b+\e))-1=deg(\b)-e-g+1+i-1=deg(\b)-e-g+i
$$
By hypotesis $d\geq 4g-2i+1$, so $2deg(\b)-e\geq 4g-2i+1$ and $deg(\b)\geq
\frac{4g-2i+1+e}{2}$. Then:
$$
\begin{array}{rl}
{h^0(\Te_X(\b+\e))-1}&{\geq \frac{4g-2i+1+e}{2}-e-g+i\geq g-\frac{e}{2}+\frac{1}{2}\geq}\\
{}&{}\\
{}&{\geq h^1(\Te_X(-\e))+\frac{1}{2}>h^1(\Te_X(-\e))}\\
\end{array}
$$
so we can apply the Lemma of Green \ref{green} and we obtain $s(\b,\b+\e)=0$.
\end{enumerate}

Now, let us see that $s(\lambda \b+\mu\e,\b+\e)=0$ when $\lambda\geq 1$ and
$0\leq
\mu\leq
\lambda-1$:
\begin{enumerate}

\item If $\lambda=1$, $\mu=0$ and it is proved.                              

\item Suppose that $\lambda>1$. Note that $(\lambda-1)\b+(\mu-1)\e$ is a nonspecial divisor, because 
 $\b$ is nonspecial and $\b+\e$ is effective.  Therefore,
$h^1(\Te_X((\lambda-1)\b+(\mu-1)\e))<h^0(\Te_X(\b+\e))-1$, because $|\b+\e|$ is base point free and
different from  $0$. Applying the Lemma of Green \ref{green} we deduce that
$s(\lambda \b+\mu\e,\b+\e)=0$.

\end{enumerate}

Finally, applying the Lemma \ref{lematecnico}, we see that:
$$
s(i\b,\b+\e,\stackrel{k-i}{\ldots},\b+\e)=s(i\b+(k-i-1)\b+\e,\b+\e)=0
$$
when $1\leq i\leq k-1$, $k\geq 2$ and the conclusion follows. \qed

\end{enumerate}

\begin{rem}\label{preciso}
If $d<4g-2i+1$ the conclusions of the Theorem could fail. Let us see an example:

Suppose that $X$ is an hyperelliptic curve of genus $g\geq 3$. Let $\b$ a very ample nonspecial divisor
of degree $2g$ and let $\b+\e=\sum_1^n g_2^1$ be a special divisor, with $0\leq n\leq g-1$.
Consider the decomposable ruled surface $\P(\Te_X(\b)\oplus \Te_X(\b+\e))$ and the unisecant linear
system $H\sim X_0+\b f$. We obtain an scroll of degree $d=2g+2m$ and speciality $i=g-n$. In this way
$d=4g-2i$. Since the scroll is decomposable, we know that:
$$
dim(s(H,H))=dim(s(\b,\b))+dim(s(\b+\e,\b+\e))+dim(s(\b+\e,\b))
$$
But $X$ is hyperelliptic and  $deg(\b)=2g$, so $(X,\Te_X(\b)$ is not projectively normal (see
\cite{green}) and $dim(s(\b,\b))>0$. Thus, the scroll and the curve of minimum
degree have not the same specialty respect to quadrics.
\end{rem}

\end{document}